\documentclass{article}
\usepackage{xypic}
\usepackage{amsthm}
\LaTeXdiagrams
\usepackage{graphicx}

\usepackage{amssymb}

\newcommand{\C}{\Bbb C}

\newcommand{\Z}{\Bbb Z}

\newtheorem{theorem}{Theorem}
\newtheorem{definition}{Definition}
\newtheorem{proposition}{Proposition}
\newtheorem{lemma}{Lemma}

\theoremstyle{remark}
\newtheorem*{remark}{Remark}

\begin{document}

\title{The Chow ring of the symmetric space $Gl(2n,\C)/SO(2n,\C)$}
\author{Rebecca E. Field }
\date{March 19, 2001}
\maketitle
\

{\footnotesize  
\begin{quote}
ABSTRACT. We show this Chow ring is 
$\Z\oplus \Z$.  
We do this by partitioning the space into $2n$ 
subvarieties
each of which is fibered over $Gl(2n-2,\C)/SO(2n-2,\C)$.
\end{quote}
}

\section{Introduction}

This paper is devoted to the following result. 

\begin{theorem} $CH^*(Gl(2n)/SO(2n)) \cong \Z \oplus \Z y,$ where y is a 
codimension $n$ cycle.
\end{theorem}
\vspace{-.2cm}
\noindent Throughout this paper $Gl(n)$, $O(n)$, $SO(n)$, etc. denote 
the complex algebraic groups of these types.

Theorem 1 should be contrasted with the following result which will used 
in its proof.

\begin{lemma} $CH^0(Gl(2n)/O(2n)) \cong \Z$ and 
$CH^i(Gl(2n)/O(2n)) = 0$ for $i \geq 1$.
\end{lemma}
\begin{proof}
$Gl(2n)/O(2n) \cong Symm_{2n}(\C)$, where $Symm_{2n}(\C)$ is the space of 
symmetric, non-degenerate $2n \times 2n$ matrices over $\C$.  
$Symm_{2n}(\C)$ 
is 
an open subset of the 
vector space 
of all $2n \times 2n$ 
symmetric 
matrices over $\C$.
By the fundamental exact 
sequence for Chow groups (Lemma 2 below), the Chow groups of any Zariski open 
subset of affine space vanish in codimensions higher than zero and are 
$\Z$ in 
codimension zero.
\end{proof}

There is much work in the literature on the geometry of symmetric
spaces of the form $G/K$ for $G$ an adjoint group and $K$ the fixed
point set of an involution of $G$, see for example \cite{DeCP}. However,
these results do not include a computation of the Chow ring. 
Moreover, the Chow ring of such a $G/K$ does not determine the Chow ring of
the symmetric space $G/K^0$ where $K^0$ is the connected component of the
identity in $K$.  This can be significantly more complicated as Theorem 1
and Lemma 1 show.

In fact, there is currently no general theorem that computes the Chow
ring of a reductive symmetric space $G/K$ in terms of Lie theoretic
data.  The most general result I know of is when $G/K$ is a group; the
Chow ring $mod \hspace{.15cm} p$ is computed 
in Kac \cite{Kac}.  For a history of
similar problems and extensive references to previous results, see the
survey article \cite{Mim}.

In a subsequent paper, we will 
use Theorem 1 as the key step in the computation of the Chow ring of the 
classifying space $BSO(2n)$.
The method of computing the Chow ring of $Gl(2n)/SO(2n)$ given in this paper 
can be extended to a computation of the cohomology of this symmetric 
space for all other cohomology theories.

\noindent{\em Acknowledgements}\\
I would like to thank Brendan Hassett, 
J. Peter May, 
Burt Totaro, 
Madhav Nori,  and especially Ian Grojnowski, without whom this paper would be 
impossible to read if it ever appeared at all.

\section{Basic results}

For notation and conventions on Chow rings, we will refer to Fulton's 
{\em Intersection theory} \cite{F} with one minor difference.  The Chow ring 
that we will be using is the ring of algebraic cycles mod rational equivalence 
(denoted $CH^*(X)$) rather than the ring of cycles mod algebraic equivalence 
(denoted $A^*(X)$) as used by Fulton.

The basic result about Chow rings that we will need is the following.

\begin{lemma} Let $Y \subseteq X$ be a closed subvariety.  Then 

$$CH_*Y \to CH_*X \to CH_*(X-Y) \to 0$$
is exact.
\end{lemma}
\vspace{-.2cm}\noindent Recall that $CH^iX=CH_{dimX-i}X.$

As in \cite{FMSS}, given an action of an 
algebraic group $\Gamma$ on a variety $X$, one forms 
$CH^\Gamma_kX = Z^\Gamma_kX/R^\Gamma_kX$.  Here $Z^\Gamma_kX$ is the free 
Abelian group generated by the $\Gamma$ stable closed subvarieties of 
$X$ and $R^\Gamma_k X$ is the 
subgroup generated by all divisors of eigenfunctions on $\Gamma$.  
A rational function $f$ on $X$ 
is an eigenfunction if $g\cdot f=\chi(g)f$ for 
all $g \in \Gamma$ and for some one dimensional character 
$\chi$ 
of $\Gamma$.  

\begin{theorem}(Fulton, MacPherson, Sotille, Sturmfels)\cite{FMSS} 
If a connected 
solvable linear algebraic group $\Gamma$ 
acts on a scheme $X$, then the canonical homorphism 
$A^\Gamma_kX \rightarrow A_kX$ is an isomorphism.  
\end{theorem}   

These results allow us to restrict attention to $B$-stable cycles in our 
proof of Theorem 1.  In fact the decomposition of $Gl(2n)/SO(2n)$ 
into $B$-orbits 
is not used essentially in the proof of the theorem, but it does 
provide a convenient notation for cycles.

\medskip
We now recall some 
results about $Gl(2n)/SO(2n)$ and its 
partition into $B$-orbits.
All of these statements are well known and their proofs are easy 
linear algebra.
We let $Symm_k(\C)$ denote the 
subset of $Gl(k)$ consisting of symmetric matrices, and define the map 
$$f:Gl(k) \rightarrow Symm_k(\C)$$ 
by 
$f(g)= gg^t$.  This map induces an isomorphism 
$$Gl(k)/O(k) \cong Symm_k(\C),$$ 
and allows us to make the identification
$$Gl(k)/SO(k) \cong 
\{(q,\epsilon  ) | q \in Gl(k)/O(k) \mbox { and } \epsilon^2 = det(q)\}.$$  

With this identification, the obvious double cover 
$$
\pi : Gl(2n)/SO(2n) \to Gl(2n)/O(2n)
$$ 
takes a pair $(q, \epsilon)$ to $q.$  We will identify 
$Symm_k(\C)$ with the 
space of symmetric non-degenerate 
bilinear forms on $\C^{2n}$ by setting $q_A(v,w)=v^tAw$, 
so that $Gl(2n)$ acts on bilinear forms by $(gq)(v,w)=q(g^tv,g^tw)$.

Define the Borel subgroup $B$ to be 
the subgroup of upper triangular matrices, i.e. 
the stabilizer of the standard flag 
$F_1 \subset F_2 \subset \cdots \subset F_{2n},$ where 
$F_i=$\linebreak$<e_1, e_2, \ldots , e_i>$.
We have an isomorphism between $B$-orbits in 
$Gl(2n)/O(2n)$ and involutions in the Weyl 
group $S_{2n}$:
$$B\backslash Gl(2n)/O(2n) \cong \{ \omega \in S_{2n} \mid \omega^2=1\},$$
defined by sending a quadratic form 
$q$ to the 
relative position between the standard flag $F$ and the flag $F^{\perp_q}$ 
of orthogonal 
complements with respect to $q$.  
It is well known and easy to see that this map is an isomorphism (this is 
a consequence of Gram-Schmidt orthonormalization).

To describe the $B$-orbits on $Gl(2n)/SO(2n)$, it is enough to describe 
the pullback of the orbits on $Gl(2n)/O(2n)$ through the double cover 
$\pi$.  An inverse image $\pi^{-1}(O)$ 
is either a single orbit or a union of two 
disjoint orbits in the double cover $Gl(2n)/SO(2n)$.

\begin{lemma}The pullback of a $B$-orbit is the union of two disjoint orbits 
if and only if the permutation indexing the orbit is fixed point free.
\end{lemma}
Let $O$ be an orbit in $Gl(2n)/SO(2n)$ and 
$q \in O$.  Then $\pi^{-1}(O)$ is a single $B$-orbit if and only if 
the stabilizer of $q$ in $B$ is disconnected; i.e. if and only if 
the stabilizer of $q$ in the set of diagonal matrices $T$ is disconnected.
Choose $q$ such that 
$q(x,y)=x^t\omega y$ where $\omega \in S_{2n}$ is the permutation 
indexing $O$.  Then $q(e_i,e_j)=\delta_{\omega(i),j}$.  The component 
group of the stabilizer of $q$ in $T$ is $(\Z/2\Z)^l$, where $l$ is the 
number of fixed points of $\omega$ on $\{1,2,\ldots,2n\}$.

\begin{definition}A $B$-orbit $O$ on $Gl(2n)/O(2n)$ is 
\textbf{fixed point free} if $\pi^{-1}(O)$ consists of two disjoint $B$-orbits.
\end{definition}

\begin{proposition} The Chow ring of $Gl(2n)/SO(2n)$ is generated by 
elements corresponding to closures of 
fixed 
point free orbits in $Gl(2n)/O(2n)$ along with 
the codimension zero cycle.
\end{proposition}

The point of this paper is to prove that all but one of these codimension 
non-zero generators 
are 
rationally equivalent to zero.

\begin{proof}
As Lemma 1
states, the Chow ring of $Gl(2n)/O(2n)$ 
is trivial.  Therefore, if $O$ is a 
$B$-orbit in $Gl(2n)/O(2n)$ whose codimension is larger than 
zero, then ${\pi}^{-1}(O)\sim 0$ in $CH^*(Gl(2n)/SO(2n))$.

By 
Definition 1,
any orbit in $Gl(2n)/O(2n)$ which is 
not fixed point free 
lifts to a single orbit in $Gl(2n)/SO(2n)$ and is therefore rationally 
equivalent to zero.

On the other hand, if an 
orbit $O$ is 
fixed point free, then 
$\pi^{-1}(O)=O_+\amalg O_-$ where 
$O_+$ and $O_-$ are two disjoint copies of $O$ distinguished by the sign of 
$\epsilon$.  
Again, since the Chow ring of the base space is trivial, 
$0 \sim \pi^{-1}(O)\sim O_+ + O_-$, so $O_+ \sim -O_-$.  Therefore, 
$CH^*(Gl(2n)/SO(2n))$ is generated by elements corresponding to closures of 
fixed 
point free orbits in $Gl(2n)/O(2n)$, along with the codimension zero cycle.
\end{proof}

\section{Proof of Theorem 1}

Define the graph of 
$B$-orbits to be a 
graph whose 
vertices are 
$B$-orbits 
and whose edges are codimension $1$ 
inclusion relations between orbit closures.  
This graph of $B$-orbits contains the full subgraph of fixed point free 
$B$-orbits.
The proof of Theorem 1 was arrived at through a careful examination of   
the fixed point free graphs of $Gl(2n)/O(2n),$ for $n \leq 4$.  
These examples 
are reproduced in the appendix to make the proof easier to visualize.

The proof is by induction.
We start by decomposing $Gl(2n)/SO(2n)$ into $2n$ disjoint
subvarieties, each of which may be compared with 
$Gl(2n-2)/SO(2n-2)$.  
Our induction hypothesis will be 
$CH^*(Gl(2n-2)/SO(2n-2)) \cong \Z x_0 \oplus \Z y$, where $x_0$ is the 
codimension zero cycle and $y$ is in codimension $n-1$.  
We will then construct a map from all but one of the $2n$ disjoint
subvarieties to $Gl(2n-2)/SO(2n-2)$.  This map will be a trivial fibration.
The remaining subvariety is easily dealt with using Proposition 1.
This map will 
take fixed point free orbits of $Gl(2n)/SO(2n)$ to fixed point free orbits 
of $Gl(2n-2)/SO(2n-2)$, 
hence the relevance of the examples in the appendix.
We will use this map and the induction step 
to show that each of the 
$2n$ subsets contributes at most $\Z \oplus \Z$ to $CH^*(Gl(2n)/SO(2n))$.  
From there, 
a second and third induction and the results of the previous section
will 
show that only two of these 
disjoint subsets actually contribute to the Chow ring.

We will start the induction with the case $n=1$.  Since there is only one 
fixed point free permutation in two letters, namely $(1 \hspace{.2cm} 2)$, 
the Chow ring $CH^*(Gl(2)/SO(2))$ has a single generator in addition to 
the trivial codimension zero cycle.  Since this single fixed point free 
orbit is the largest fixed point free orbit, it cannot be the zero of any 
$B$-semi-invariant function.  Therefore, 
$CH^*(Gl(2)/SO(2)) \cong \Z x_0 \oplus \Z y$, where $x_0$ is the codimension 
zero cycle and $y$ is this codimension $1$ cycle.

We now define a decomposition of $Gl(2n)/SO(2n)$ into 
$2n$ disjoint subspaces.  

Define $X_i$ to be the subvariety of $Gl(2n)/SO(2n)$ consisting of pairs 
$(q,\epsilon)$ where 
$q(e_{2n},e_{2n})=q(e_{2n-1},e_{2n})=\cdots=q(e_{i+1},e_{2n})=0$ and 
$q(e_i,e_{2n})\neq 0$.
The following properties of the varieties $X_i$ are immediate.

\begin{lemma}
$$Gl(2n)/SO(2n)=X_1 \amalg X_2 \amalg \cdots 
\amalg X_{2n}.$$  
Also,
$$\overline{X_i}=\coprod_{j\geq i}X_j.$$
The $X_i$ and the $\overline{X_i}$ are all quasi-projective varieties, and 
each $X_i$ is $B$-stable.
\end{lemma}

Notice that if $(q, \epsilon) \in X_{2n}$, then the $B$-orbit through $q$ is 
not fixed point free, and hence by Proposition 1, the only $B$-orbit in 
$X_{2n}$ that contributes to the Chow ring is the 
codimension zero 
open orbit.

For
$i<2n$, we define a map 
$$f_i:X_i \longrightarrow Gl(2n-2)/SO(2n-2)$$
as follows.
For $(q,\epsilon)\in X_i$, 
we know that 
$q(e_{i+1},e_{2n})=q(e_{i+2},e_{2n})=\cdots=q(e_{2n},e_{2n})=0$, while 
$q(e_i,e_{2n})\neq 0$.  Therefore, the quadratic form $q$ is nondegenerate 
on $<e_i,e_{2n}>$, the subspace of $\C^{2n}$ generated by $e_i$ and $e_{2n}$.  
The orthogonal compliment $<e_i,e_{2n}>^\perp$ 
with respect to $q$ of this subspace is 
isomorphic to $\C^{2n-2}$. 

Let $q'$ be the quadratic form on $\C^{2n-2}$ defined by restricting 
$q$ to $<e_i,e_{2n}>^\perp$ and let 
$\epsilon'=\epsilon/q(e_i,e_{2n})\sqrt{-1}$.  As 
$det(q)=-q(e_i,e_{2n})^2det(q')$, we know $(\epsilon')^2=det(q').$
Define $f_i(q,\epsilon)=(q',\epsilon').$

This map takes $B$-orbits to $\overline{B}$-orbits, where $\overline{B}$ is 
the stabilizer in $Gl(2n-2)$ of the standard flag associated to the basis
$$\overline{e_1},\ldots,\overline{e_{i-1}},\overline{e_{i+1}},\ldots,\overline{e_{2n-1}},$$ 
where 
$\overline{e_j}$ is equal to the orthogonal projection of 
$e_j$ onto $<e_i,e_{2n}>^\perp$.

\begin{lemma}$f_i$ is a trivial fibration with fibers  
isomorphic to $\C^*\times \C^{2n+i-2}.$
\end{lemma}
\begin{proof}
Fix a quadratic form $q'$
in $Gl(2n-2)/SO(2n-2)$.  An element 
$q \in f^{-1}(q')$
is determined by the numbers 
$${q(e_j,e_{2n})\mbox{ for }j<i, 
\mbox{ and } \atop q(e_j,e_i)\mbox{ for }j\neq i.}$$
These numbers may be chosen freely subject only to the constraint that 
$q(e_i,e_{2n})\neq 0.$ 
\end{proof}

\begin{lemma}$CH^*(X_i)$ is a quotient of $\Z \oplus \Z$.
\end{lemma}
\begin{proof}

Since the map $f_i$ is a fibration with fibers isomorphic to 
an open subset of affine space, it induces 
a surjection 
of 
Chow rings (\cite{F} example 1.9.2).  
By our induction hypothesis 
$CH^*(Gl(2n-2)/SO(2n-2))\cong \Z \oplus \Z$.

More specifically, the Chow ring of $Gl(2n-2)/SO(2n-2)$ is generated by 
cycles $x_0$ and $y$,
where $x_0$ is the codimension zero cycle and 
$y$ is in codimension $n-1$.  
Therefore, $CH^*(X_i)$ is generated by the pullbacks 
of these cycles.
\end{proof}

\begin{remark}
In fact $CH^*(X_i)\cong \Z \oplus \Z$ as follows from the proof of Theorem 1 
below.  This can also easily be seen directly.
\end{remark}

\begin{proposition} If
$n>1$, then 
$CH^*(Gl(2n)/SO(2n))=CH^*(\overline{X_{2n}})$ is 
a quotient of 
$\Z \widetilde{x_0} \oplus \Z f_{2n-1}^*(y) \oplus \Z f_{2n-2}^*(y) 
\oplus \cdots \oplus \Z f_1^*(y)$, where $\widetilde{x_0}$ is the codimension 
zero cycle in $\overline{X_{2n}}$. 
\end{proposition} 

\begin{proof}
We will show by induction on $i$ that $CH^*(\overline{X_i})$ is 
a quotient of 
$\Z f_{i}^*(y) \oplus \Z f_{i-1}^*(y) 
\oplus \cdots \oplus \Z f_1^*(y).$

For the case $i=1$, since $\overline{X_1}=X_1$, Lemma 6 tells us that   
$CH^*(\overline{X_1})$ is a quotient of $\Z f^*_1(x_0) \oplus \Z f^*_1(y)$.  
The pullback of the codimension zero orbit $f^*_1(x_0)$ is indexed by the 
permutation $(1 \hspace{.2cm} 2n)$ so since $n>1$ it 
is not a fixed point free orbit 
in $Gl(2n)/SO(2n)$.  Therefore, by Proposition 1, $f_1^*(x_0)$ 
is rationally equivalent to zero in 
$CH^*(Gl(2n)/SO(2n))$, and the subvariety 
$\overline{X_1}$ contributes at most $f^*_1(y)$
to the Chow ring of the whole symmetric space.  This completes the $i=1$ case.

As stated in 
Lemma 2, we know that 
$$CH^*(\overline{X_{i-1}}) \to CH^*(\overline{X_i}) \to 
CH^*(X_i) \to 0$$
is exact.  Therefore, by induction, $CH^*(\overline{X_i})$ is 
a quotient of 
$\Z f_i^*(x_0) \oplus \Z f_i^*(y) \oplus \Z f_{i-1}^*(y) \oplus \cdots 
\oplus \Z f_1^*(y)$.  Again, we note that 
$f_i^*(x_0)$ is indexed by the permutation 
$(i \hspace{.2cm} 2n)$ so is not fixed point free.  For $i<2n$, this orbit 
has codimension larger than zero, so by Proposition 1, this orbit is 
rationally equivalent to zero in $CH^*(Gl(2n)/SO(2n))$.  
Finally, if $i=2n$, as we noticed previously (immediately following Lemma 4), 
the subvariety $X_{2n}$ contributes only the 
codimension zero cycle $\widetilde{x_0}$.  
\end{proof}

Note that we have reduced the generators down to the largest fixed point free 
$B$-orbit in each of the $X_i$.  These $2n-1$ $B$-orbits are visible 
in the appendix along the bottom right hand side of each graph.

We claim that $f_j^*(y)$ is rationally equivalent to zero in 
$CH^*(Gl(2n)/SO(2n))$ for $j<2n$.
To this end, we 
inductively show that $CH^*(\overline{X_j})$ is 
a quotient of
$\Z  f_j^*(y)$ for $j < 2n$.

Define the function $g_j:\overline{X_j} \longrightarrow \C$ 
by $g_j(q)=q(e_j,e_{2n})$.  
This function is clearly non-zero on 
$X_j$ and is clearly zero on 
all of $\overline{X_{j-1}}$. 
To see that this zero is a simple one, we look at the function 
restricted to a line in $Gl(2n)/SO(2n)$.  
Let $q \in X_j, q'\in X_{j-1}$ and consider the line of quadratic forms 
$\{q'+aq\mid a \in \C\}$.
Clearly the function $g_j$ 
has a 
simple zero along this line. 
This line meets $X_{j-1}$
only when 
$a=0$ and is otherwise completely 
contained in $X_j$.  
It is clear that the line 
intersects $\overline{X_{j-1}}$ transversally since $\overline{X_{j-1}}$ 
is an open subvariety of the vector 
space of quadratic forms for which $q(e_j,e_{2n})=0$.  
Therefore, the intersection is transversal and the zero is simple.

By our induction step, $CH^*(\overline{X_{j-1}})$ is 
a quotient of 
$\Z f_{j-1}^*(y)$, so $CH^*(\overline{X_j})$ is a quotient of 
$\Z f_{j-1}^*(y)\oplus \Z f_{j}^*(y)$.
Let $\gamma$ be the $B$-orbit in $Gl(2n-2)/SO(2n-2)$ 
representing the cycle $y$.  Restricting $g_j$ to 
the closure of 
$f_j^{-1}(\gamma)$ in $\overline{X_j}$, this argument shows that 
$f_{j-1}^*(y) \sim 0$ in $CH^*(\overline{X_j})$.

We have shown $f_j^*(y)\sim 0$ for $j<2n-1$.  Consider $f_{2n-1}^*(y)$.
This cycle is the 
closure of the 
largest fixed point free orbit, so
any $B$-semi-invariant function vanishing along 
it
must be defined 
on the closure of a non-fixed point free orbit.  Therefore, any 
such function will produce a 
relation involving both copies of the orbit pulled back from $Gl(2n)/O(2n)$.  
Since we already know that the sum of these copies is rationally equivalent 
to zero, such a function can produce no new relations.  The same argument 
shows that no multiple of $f_{2n-1}^*(y)$ is rationally equivalent to zero.
 
We have shown that $f_{j}^*(y) \sim 0$ for $1<j<2n-1$ and that 
no multiple of $f_{2n-1}^*(y)$ is rationally equivalent to zero.  
This combined with  
Proposition 2 proves Theorem 1.

\section{Appendix}

\hspace{.4cm}  As mentioned in the proof, 
the graph of fixed point free $B$-orbits for \linebreak
$Gl(2)/SO(2)$ is a single point corresponding to the orbit $O_{(12)}$.

The graph of fixed point free $B$-orbits for $Gl(4)/SO(4)$ is:
$$\spreaddiagramrows{-.3pc}
\xymatrix{
O_{(14)(23)} \ar[d]\\
O_{(13)(24)} \ar[d]\\
O_{(12)(34)}
}
$$ 

\noindent where the codimensions of the orbits are (counting from the 
bottom) 2, 3, and 4 and the arrows represent inclusion of an orbit 
in the closure of the larger orbit.  Since the graph for $Gl(2)/SO(2)$ is 
a single point, this is also its decomposition into subgraphs.

This leads to the induction step for the graph of fixed point free 
$B$-orbits for $Gl(6)/SO(6)$, namely: 
\vspace{-.2cm}
$$\spreaddiagramcolumns{-1pc}
\spreaddiagramrows{-.3pc}
\xymatrix{
&&O_{(16)(25)(34)} \morphism\dashed\tip\notip[1,-1] \ar[dr]&&\\
&O_{(15)(26)(34)}\morphism\dashed\tip\notip[1,-1]\ar[dr]&&O_{(16)(24)(35)}\morphism\dashed\tip\notip[1,-3]\morphism\dashed\tip\notip[1,-1]\ar[dr]&\\ 
O_{(15)(24)(36)}\morphism\dashed\tip\notip[1,0]\ar[drr]&&O_{(14)(26)(35)}\morphism\dashed\tip\notip[1,0]\ar[drr]&&O_{(16)(23)(45)}\morphism\dashed\tip\notip[1,-4]\ar[d]\\
O_{(15)(23)(46)}\morphism\dashed\tip\notip[1,0]\ar[drr]&&O_{(14)(25)(36)}\morphism\dashed\tip\notip[1,0]\morphism\dashed\tip\notip[1,-2]\ar[drr]&&
O_{(13)(26)(45)}\morphism\dashed\tip\notip[1,-2]\ar[d]\\
O_{(14)(23)(56)}\ar[dr]&&O_{(13)(25)(46)}\morphism\dashed\tip\notip[1,-1]\ar[dr]&&O_{(12)(36)(45)}\ar[dl]\\
&O_{(13)(24)(56)}\ar[dr]&&O_{(12)(35)(46)}\ar[dl]&\\
&&O_{(12)(34)(56)}&&
}
\vspace{-.1cm}
$$

\noindent where the solid diagonal 
lines that are more or less parallel to each other 
represent inclusion relations within the 
$X_i$ and the dotted lines represent other inclusion relations that 
are not relevant to 
our proof
Furthermore, orbits that appear 
on the same row have the same codimension.

The following is the graph of
fixed point free $B$-orbits in 
$Gl(8)/O(8)$.  Notice its seven subgraphs. 

\nopagebreak
\hspace{-.2in}
\includegraphics[height =3.2in]{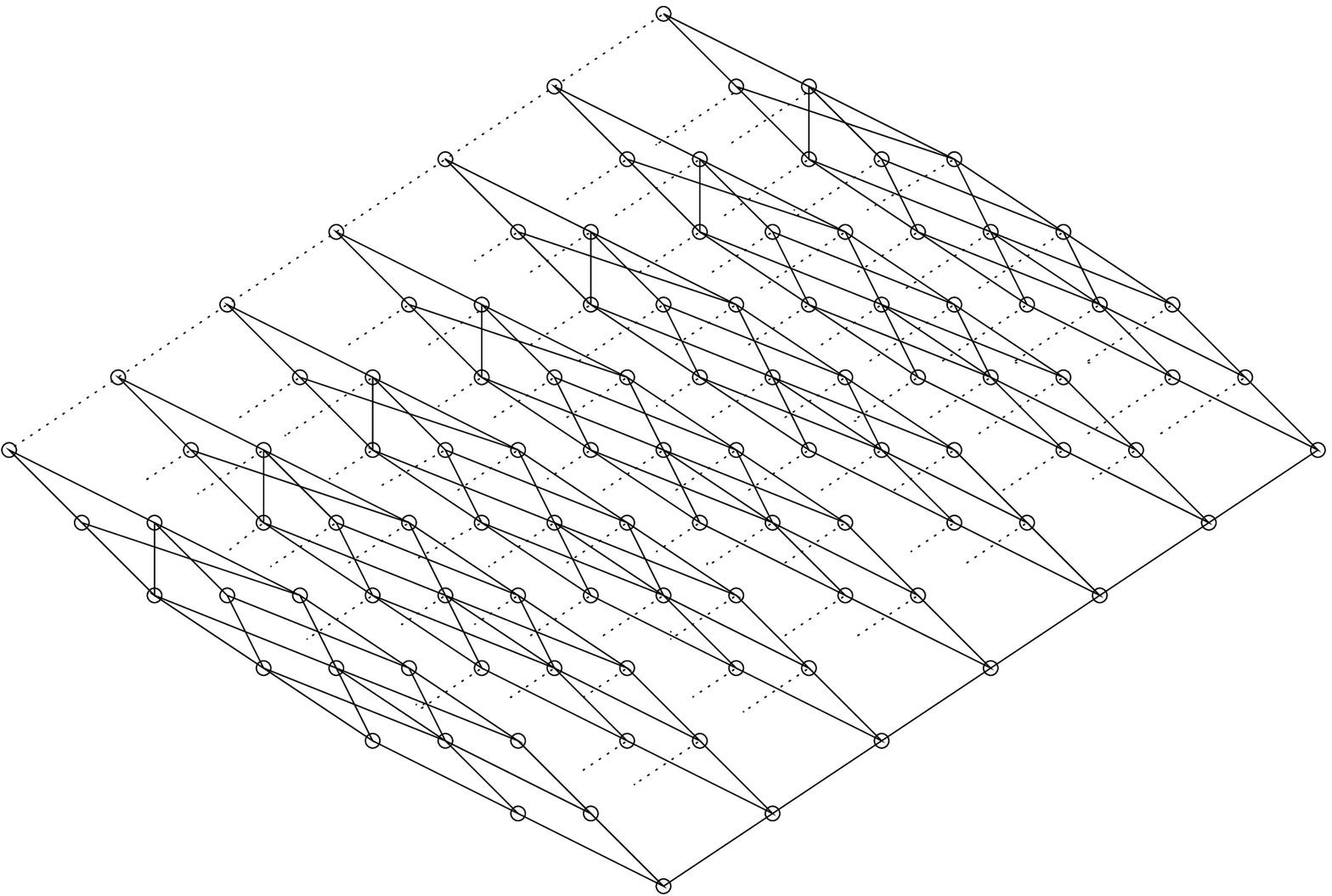}

{\scriptsize  
DEPARTMENT OF \vspace{-.13cm} MATHEMATICS, UNIVERSITY OF WISCONSIN, 480 LINCOLN DRIVE, 
\vspace{-.075cm}
MADISON, WI 53706-1388 USA}


{\small {\em E-mail address:}
\verb|field@math.wisc.edu|}

\end{document}